\magnification=\magstep1
\tolerance=2000

\baselineskip=17pt 

\font\tenopen=msbm10
\font\sevenopen=msbm7
\font\fiveopen=msbm5
\newfam\openfam
\def\openo{\fam\openfam\tenopen}
\textfont\openfam=\tenopen
\scriptfont\openfam=\sevenopen
\scriptscriptfont\openfam=\fiveopen

\font\tenujsym=msam10
\font\sevenujsym=msam7
\font\fiveujsym=msam5
\newfam\ujsymfam

\textfont\ujsymfam=\tenujsym
\scriptfont\ujsymfam=\sevenujsym
\scriptscriptfont\ujsymfam=\fiveujsym

\def\C{{\openo C}}
\def\R{{\openo R}}
\def\N{{\openo N}}

\def\loc{{\rm loc}}
\def\e{\varepsilon}

\centerline{\bf Necessary and sufficient Tauberian conditions for} 

\centerline{\bf the logarithmic summability of functions and sequences} 

\bigskip

\centerline{By Ferenc M\'oricz} 

\bigskip

{\it University of Szeged, Bolyai Institute, Aradi v\'ertan\'uk tere} 1, 6720 {\it Szeged, Hungary}

e-mail: moricz@math.u-szeged.hu
\bigskip

\centerline{\it Abstract} 

\bigskip

\noindent Let $s: [1, \infty) \to \C$ be a locally integrable function in Lebesgue's sense on the 
infinite interval $[1, \infty)$. 
We say that $s$ is summable $(L, 1)$ if there exists some $A\in \C$ such that 
$$\lim_{t\to \infty} \tau(t) = A, \quad {\rm where} \quad \tau(t):= {1\over \log t} 
\int^t_1 {s(u) \over u} du.\leqno(*)$$
It is clear that if the ordinary limit $s(t) \to A$ exists, then the limit 
$\tau(t) \to A$ also exists as $t\to \infty$. 
We present sufficient conditions, which are also necessary in order that the converse implication 
hold true. As corollaries, we obtain  so-called Tauberian theorems which are analogous to 
those known in the case of summability $(C,1)$. 
For example, if the function $s$ is slowly oscillating, by which we mean that for every $\e>0$ there exist 
$t_0 = t_0 (\e) > 1$ and $\lambda=\lambda(\e) > 1$ such that 
$$|s(u) - s(t)| \le \e \quad {\rm whenever}\quad t_0 \le t < u \le t^\lambda,$$
then the converse implication holds true: the ordinary 
convergence $\lim_{t\to \infty} s(t) = A$ follows from ($*$). 

We also present necessary and sufficient Tauberian conditions under which 
the ordinary convergence of a numerical sequence $(s_k)$ follows from its logarithmic summability. Among 
others, we give a more transparent  proof of an earlier Tauberian theorem due to Kwee [3]. 

\vfill\eject

\centerline{1. {\it Introduction: Summability} $(C, 1)$ {\it and} $(L, 1)$ {\it of functions} }

Let $s:[0, \infty) \to \C$ be an integrable function in Lebesgue's sense on every bounded interval 
$[0,t]$, $t>0$, in symbols: $s\in L_{\loc} [0, \infty)$. We recall (see, e.g., [2, p. 11]) that the function 
$s$ is said to be Ces\`aro summable of first order, briefly: summable $(C, 1)$, if there exists some 
$A\in \C$ such that 
$$\lim_{t\to \infty} \sigma(t) = A, \quad {\rm where} \quad \sigma(t):= {1\over t} \int^t_0 s(u) du, 
\quad t>0. \leqno(1.1)$$

Clearly, if the ordinary limit 
$$\lim_{t\to \infty} s(t) = A\leqno(1.2)$$
exists, then the limit in (1.1) also exists. The converse 
implication holds true only under some supplementary, so-called Tauberian condition(s). 

We note that the left endpoint of the definition domain of the function $s$ is indifferent in (1.1). 
That is, given any $a>0$, the existence of the limit in (1.1) is equivalent with the existence of the following 
one: 
$$\lim_{t\to \infty} {1\over t} \int^t_a s(u) du = A.$$

As to the Ces\`aro summability of order $\alpha$, where $\alpha\ge 0$ is a real number, briefly: 
summability $(C, \alpha)$, we refer to [7, p. 26]. The case $\alpha=0$ is ordinary convergence. 

Next, let $s: [1, \infty) \to \C$ be such that $s\in L_{\loc} [1, \infty)$. Motivated by the 
concept of logarithmic (sometimes also called harmonic) summability of a numerical 
sequence (see, e.g., in [5]), the function $s$ is said to be logarithmic summable of first order, 
briefly: summable $(L, 1)$, if there exists some $A\in \C$ such that 
$$\lim_{t\to \infty} \tau(t) = A, \quad {\rm where}\quad \tau(t):= {1\over \log t} 
\int^t_1 {s(u)\over u} du, \quad t>1,\leqno(1.3)$$
where the logarithm is to the natural base $e$. 

In Section 4, we will prove that summability $(C, 1)$ of a function implies its 
summability $(L, 1)$ to the same limit, but the converse implication is not true in general. 

We note that a complex-valued function $s\in L_{\loc} [e, \infty)$ is said to be logarithmic summable of order 2, briefly: 
summable $(L, 2)$, if there exists some $A\in \C$ such that 
$$\lim_{t\to \infty} \tau_2(t) = A\quad {\rm where}\quad 
\tau_2(t):= {1\over \log \log t} \int^t_e {s(u)\over u\log u} du, \quad t>e.\leqno(1.4)$$

We also note that in the particular cases when 
$$s(u):= \int^u_0 f(x) dx, \ u>0; \quad {\rm or}\quad s(u) := \int^u_1 f(x) dx, \quad u>1;\leqno(1.5)$$
where $f$ is a locally integrable function on $[0, \infty)$ or $[1, \infty)$, respectively, 
the above summability methods may be applied to assign value to the following integrals, respectively: 
$$\int^\infty_0 f(x) dx \quad {\rm or} \quad \int^\infty_1 f(x) dx.\leqno(1.6)$$

If the finite limit in (1.2) exists, then the improper integrals 
$\int^{\to \infty}_0 f(x) dx$ and $\int^{\to \infty}_1 f(x) dx$ exist, respectively. In the case when  the finite limit in (1.1) 
exists, then the integral in (1.6) (i) is said to be summable $(C,1)$; while in the case when only the finite limit in 
(1.3) exists, then the integral in (1.6) (ii) is said to be summable $(L, 1)$. 

\vglue1cm

\centerline{2. {\it Main results} }

In our first new result we characterize the converse implication when the ordinary limit of a real-valued 
function at $\infty$ follows from its summability $(L, 1)$. 

\bigskip

THEOREM 1. {\it If a real-valued function $s\in L_{\loc} [1, \infty)$ is summable $(L, 1)$ 
to some $A\in \R$, then the ordinary limit} (1.2) {\it exists if and only if 
$$\limsup_{\lambda\to 1+} \ \liminf_{t\to \infty} {1\over (\lambda-1)\log t} \int^{t^\lambda}_t 
{s(u) - s(t)\over u} du \ge 0\leqno(2.1)$$
and}
$$\limsup_{\lambda\to 1-} \ \liminf_{t\to \infty} {1\over (1-\lambda) \log t} \int^t_{t^\lambda} {s(t) - s(u)\over u} 
du \ge 0.\leqno(2.2)$$

Motivated by the definition of the `slow decrease' with respect to summability $(C, 1)$ (see, e.g., [2, pp. 124-125; 
and cf. our Remark 1 below]), 
we say that a function $s: [1, \infty) \to \R$ is slowly decreasing with respect to summability $(L, 1)$ if 
for every $\e>0$ there exist $t_0=t_0(\e) > 1$ and $\lambda=\lambda(\e)>1$ such that 
$$s(u)-s(t) \ge -\e\quad {\rm whenever} \quad t_0 \le t<u\le t^\lambda.\leqno(2.3)$$

It is easy to check that a function $s$ is slowly decreasing with respect to summability $(L, 1)$ 
if and only if 
$$\lim_{\lambda\to 1+} \ \liminf_{t\to \infty} \ \inf_{t< u\le t^\lambda} (s(u) - s(t)) \ge 0.\leqno(2.4)$$
Since the auxiliary function 
$$a(\lambda): = \liminf_{t\to \infty} \ \inf_{t<u\le t^\lambda} (s(u) - s(t))$$
is evidently decreasing in $\lambda$ on the infinite interval $(1, \infty)$, the right limit in (2.4) exists and 
$\lim_{\lambda\to 1+}$ can be replaced by $\sup_{\lambda>1}$. 

It is clear that if a function $s\in L_{\loc} [1, \infty)$ is slowly decreasing with respect to summability 
$(L,1)$, then conditions (2.1) and (2.2) are trivially satisfied. Thus, the next corollary is an immediate consequence of 
Theorem 1. 

\bigskip

COROLLARY 1. {\it Suppose a real-valued function $s\in L_{\loc} [1, \infty)$ is slowly decreasing 
with respect to summability $(L, 1)$. If $s$ is summable $(L, 1)$ to some $A\subset \R$, then the ordinary 
limit} (1.2) {\it also exists.} 

Historically, the term `slow decrease' was introduced by Schmidt [6] in the case of the 
summability $(C, 1)$ of sequences of real numbers. 

In our second new result we characterize the converse implication when the ordinary convergence of a 
complex-valued function follows from its summability $(L, 1)$. 

\bigskip

THEOREM 2. {\it If a complex-valued function $s\in L_{\loc} [1, \infty)$ is summable $(L, 1)$ to some 
$A\in \C$, then the ordinary limit} (1.2) {\it exists if and only if }
$$\liminf_{\lambda\to 1+} \ \limsup_{t\to \infty} \Big|{1\over (\lambda-1)\log t} \int^{t^\lambda}_t {s(u) - s(t)\over u} 
du\Big| = 0.\leqno(2.5)$$

Motivated by the definition of the `slow oscillation' with respect to summability $(C, 1)$ of numerical sequences 
introduced by Hardy [1] (see also in [2, pp. 124-125]), we say that a function $s: [1, \infty) \to \C$ is slowly oscillating with 
respect to summability $(L, 1)$ if for every $\e>0$ there exist 
$t_0 = t_0 (\e)>1$ and $\lambda=\lambda(\e)>1$ such that 
$$|s(u) - s(t)| \le \e \quad {\rm whenever}\quad t_0 \le t < u \le t^\lambda.\leqno(2.6)$$

It is easy to check that a function $s$ is slowly oscillating with respect to summability $(L, 1)$ 
if and only if 
$$\lim_{\lambda\to 1+} \limsup_{t\to \infty} \ \sup_{t< u \le t^\lambda} |s(u) - s(t)| = 0.\leqno(2.7)$$

It is clear that if a function $s\in L_{\loc} [1, \infty)$ is slowly oscillating with respect to summability $(L, 1)$, then condition 
(2.5) is trivially  satisfied. Thus, the next corollary is an immediate consequence of Theorem 2. 

\bigskip

COROLLARY 2. {\it Suppose a complex-valued function $s\in L_{\loc} [1, \infty)$ is slowly oscillating with respect to summability 
$(L, 1)$. If $s$ is summable $(L, 1)$ to some $A\in \C$. then the ordinary limit} (1.2) {\it also exists.} 

\bigskip

REMARK 1. According to Hardy's definition (see [2, pp. 124-125]), a function $s: (0, \infty) \to \C$ is said to be slowly 
oscillating if 
$$\lim(s(u) - s(t)) = 0\quad {\rm whenever}\quad u> t\to \infty\quad {\rm and} \quad u/t\to 1;\leqno(2.8)$$
and a function $s: (0, \infty) \to \R$ is said to be slowly decreasing if 
$$\liminf(s(u) - s(t)) \ge 0 \quad {\rm under\ the \ same\ circumstances.}\leqno(2.9)$$

We claim that definition (2.8) is equivalent to the following one: for every $\e>0$ there exist 
$t_0 = t_0 (\e) > 0$ and $\lambda = \lambda(\e)>1$ such that 
$$|s(u) - s(t)|\le \e \quad {\rm whenever}\quad t_0 \le t < u \le \lambda t.\leqno(2.10)$$

The implication (2.8) $\Rightarrow$ (2.10) is trivial. To justify the converse implication (2.10) $\Rightarrow$ (2.8), 
let $\lambda>1$ be arbitrarily close to 1 and set $\e: = \log \lambda$. Then by (2.10), we have 
$$|s(u) - s(t)| \le \e \quad {\rm whenever}\quad u> t\ge t_0\quad {\rm and}\quad 
0<\log {u\over t} \le \log \lambda=\e.$$
Now, the equivalence of the two definitions claimed above is obvious. 

It is worth to consider the special case (1.6) (ii), where $f\in L_{\loc} [1, \infty)$. If $f$ is a 
real-valued function and 
$$x (\log x) f(x) \ge - C \quad {\rm at\ almost\ every} \quad x>x_0,\leqno(2.11)$$
where $C>0$ and $x_0 \ge 1$ are constants, then $s$ defined in (1.5) (ii)  is slowly decreasing with respect to summability 
$(L, 1)$, and Corollary 1 applies. 
Likewise, if $f$ is a complex-valued function and 
$$x(\log x)| f(x)| \le C \quad {\rm at\ almost\ every}\quad x>x_0,\leqno(2.12)$$
where $C>0$ and $x_0 \ge 1$ are constants, then $s$ is slowly oscillating with respect to summability 
$(L, 1)$, and Corollary 2 applies. 

Condition (2.11) is called a one-sided Tauberian condition, while (2.12) is called a two-sided 
Tauberian condition with respect to summability $(L,1)$. These terms go back to Landau [4] with respect to summability 
$(C,1)$ of sequences of real numbers, as well as to Hardy [1] (see also [2, p. 149]) with respect to summability 
$(C,1)$ of sequences of complex numbers. 

We note that such theorems containing appropriate additional conditions such as (2.11), (2.12), etc. are called 
`Tauberian', after A. Tauber, who first proved one of the simplest of this kind; and these
 supplementary conditions are 
called  `Tauberian conditions'. 

\vglue1cm
\centerline{3. {\it Proofs of Theorems} 1 {\it and} 2} 

The following two representations of the difference $s(t) - \tau(t)$ will be of vital importance in 
our proofs below. 

LEMMA 1. (i) {\it If $\lambda>1$ and $t>1$, then}
$$s(t) - \tau(t) = {\lambda\over \lambda-1} (\tau(t^\lambda) - \tau(t)) - {1\over (\lambda-1)\log t} 
\int^{t^\lambda}_t {s(u) - s(t)\over u} du.\leqno(3.1)$$

(ii) {\it If $0<\lambda < 1$ and $t>1$, then }
$$s(t) - \tau(t) = {\lambda\over 1-\lambda} (\tau(t) - \tau(t^\lambda)) + {1\over (1-\lambda) \log t} 
\int^t_{t^\lambda} {s(t) - s(u)\over u} du.\leqno(3.2)$$

{\it Proof. Part} (i). By definition in (1.3), we have 
$$\tau(t^\lambda) -\tau(t) = {1\over \lambda\log t} \int^{t^\lambda}_1 {s(u)\over u} du - {1\over \log t} 
\int^t_1 {s(u)\over u} du$$
$$={1-\lambda\over \lambda\log t} \int^t_1 {s(u)\over u} du + {1\over \lambda\log t} \int^{t^\lambda}_t 
{s(u)\over u} du$$
$$=-{\lambda-1\over \lambda} \tau(t) + {1\over \lambda \log t} \int^{t^\lambda}_t {s(u)\over u} du.$$
Multiplying both sides by $\lambda/(\lambda-1)$ gives 
$${\lambda\over \lambda-1} (\tau(t^\lambda) - \tau(t)) = - \tau(t) + 
{1\over (\lambda-1)\log t} \int^{t^\lambda}_t {s(u)\over u} du$$
$$=s(t) - \tau(t) + {1\over (\lambda-1) \log t} \int^{t^\lambda}_t {s(u) - s(t)\over u} du,$$
whence (3.1) follows. 

{\it Part} (ii). The proof of (3.2) is analogous to that of (3.1). 

{\it Proof of Theorem} 1. {\it Necessity.} Suppose that (1.2) is satisfied. By (1.2) and (1.3), we have 
$$\lim_{t\to \infty} (s(t) - \tau(t)) = 0\quad {\rm and} \quad 
\lim_{t\to \infty} (\tau(t^\lambda) - \tau(t)) = 0\leqno(3.3)$$
for each fixed $\lambda>1$. By (3.1) and (3.3), we conclude that 
$$\lim_{t\to \infty} {1\over \log t} \int^{t^\lambda}_t {s(u) - s(t)\over u} du =0\leqno(3.4)$$
for every $\lambda>1$. This proves (2.1) even in a stronger form. 

An analogous argument yields (2.2) for every $0<\lambda<1$ also in a stronger form. 

{\it Sufficiency}. Suppose that (2.1) and (2.2) are satisfied.
By (2.1), there exists a sequence $\lambda_j\downarrow 1$ such that 
$$\lim_{j\to \infty} \ \liminf_{t\to \infty} \ {1\over (\lambda_j -1)\log t} \int^{t^{\lambda_j}}_t 
{s(u) - s(t)\over u} du \ge 0.\leqno(3.5)$$
By (1.3), (3.1) and (3.5), we conclude that 
$$\limsup_{t\to \infty} (s(t) - \tau(t)) \le \lim_{j\to \infty} 
\limsup_{t\to \infty} {\lambda_j \over \lambda_j - 1} (\tau(t^{\lambda_j}) - \tau(t))\leqno(3.6)$$
$$+\lim_{j\to \infty} \ \limsup_{t\to \infty} \Big(-{1\over (\lambda_j -1) \log t} 
\int^{t^{\lambda_j}}_t {s(u) - s(t)\over u} du\Big)$$
$$=-\lim_{j\to \infty} \ \liminf_{t\to \infty} {1\over (\lambda_j - 1) \log t} \int^{t^{\lambda_j}}_t 
{s(u) - s(t)\over u} du \le 0.$$

By (2.2), there exists a sequence $0<\lambda_k\uparrow 1$ such that 
$$\lim_{k\to \infty} \ \liminf_{t\to \infty} {1\over (1-\lambda_k) \log t} \int^t_{t^{\lambda_k}} 
{s(t) - s(u)\over u} du \ge 0.\leqno(3.7)$$
By (1.3), (3.2) and (3.7), we conclude that 
$$\liminf_{t\to \infty} (s(t) - \tau(t)) \ge \lim_{k\to \infty} \ \liminf_{t\to \infty} 
{\lambda_k\over 1-\lambda_k}  
(\tau(t) - (\tau (t^{\lambda_k}))\leqno(3.8)$$
$$+\lim_{k\to \infty} \liminf_{t\to \infty} {1\over (1-\lambda_k) \log t} 
\int^t_{t^{\lambda_k}} {s(t)-s(u)\over u} du$$
$$=\lim_{k\to \infty} \ \liminf_{t\to \infty} {1\over (1-\lambda_k) \log t} \int^t_{t^{\lambda_k}} 
{s(t) - s(u)\over u} du \ge 0.$$

Combining (3.6) and (3.8) yields (3.3) (i), and a fortiori, we get (1.2) to be proved, due to summability 
$(L, 1)$ of the function $s$. 

The proof of Theorem 1 is complete. 

\bigskip

{\it Proof of Theorem} 2. It also hinges on Lemma 1 and runs along similar lines to the proof of 
Theorem 1. The details are left to the reader. 

\vglue1cm

\centerline{4. {\it Inclusions} }

We will prove that summability $(L,1)$ is more effective than summability $(C,1)$. 

\bigskip

THEOREM 3. {\it If a complex-valued function $s\in L_{\loc} [1, \infty)$ is summable $(C,1)$ to some 
$A\in \C$, then $s$ is also summable $(L,1)$ to the same A. The converse implication is not true in general. 

\bigskip

Proof.} (i) First, let $t:= m$, where $m=2,3,\ldots$. 
By definition in (1.4) and applying the Second Mean-Value Theorem, we get 
$$\tau(m) \log m = \sum^{m-1}_{k=1} \int^{k+1}_k {s(u)\over u} du\leqno(4.1)$$
$$=\sum^{m-1}_{k=1} \Big({1\over k} \int^{\xi_k}_k s(u) du + {1\over k+1} \int^{k+1}_{\xi_k} s(u) du\Big)$$
$$=\int^{\xi_1}_1 s(u) du + \sum^{m-2}_{k=1} {1\over k+1} \int^{\xi_{k+1}}_{\xi_k} 
s(u) du + {1\over m} \int^m_{\xi_{m-1}} s(u) du,$$
where
$$k<\xi_k<k+1\quad {\rm for}\quad k=1,2,\ldots, m-1.\leqno(4.2)$$

By definition in (1.1), we may write that 
$$\int^\eta_\xi s(u) du = \eta\sigma(\eta)-\xi \sigma(\xi), \quad 0<\xi<\eta.\leqno(4.3)$$
Making use of this equality, from (4.1) it follows that 
$$\tau(m) \log m = -(\xi_1 \sigma(\xi_1)-\sigma(1))$$
$$+\sum^{m-2}_{k=1} {1\over k+1} (\xi_{k+1} \sigma(\xi_{k+1}) - \xi_k \sigma(\xi_k)) + {1\over m} (m \sigma(m) - \xi_{m-1} \sigma(\xi_{m-1}))$$
$$=-\sigma(t) + \sum^{m-1}_{k=1} {1\over k(k+1)} \xi_k \sigma(\xi_k) + \sigma(m) - \sigma(1),$$
whence we get 
$$\tau(m) = {1\over \log m} \sum^{m-1}_{k=1} {\xi_k\over k(k+1)} \sigma(\xi_k) + 
{1\over \log m} (\sigma(m) - \sigma(1)).\leqno(4.4)$$

We will apply Toeplitz' theorem on the summability of numerical sequences (see, e.g., [8, p. 74]) in the 
case of (4.4) with the  infinite triangular matrix 
$$\Big(a_{m,k}:= {1\over \log m} {\xi_k\over k(k+1)}, \ k=1,2,\ldots, m-1; m=2,3,\ldots\Big).$$
By (4.2), we have 
$${1\over \log m} \sum^{m-1}_{k=1} {1\over k+1} < \sum^{m-1}_{k=1} a_{m,k} < {1\over \log m} 
\sum^{m-1}_{k=1} {1\over k}, \quad  m=2,3,\ldots;$$
whence it follows that 
$$\lim_{m\to \infty} \ \sum^{m-1}_{k=1} a_{m,k} =1.$$
It is also clear that 
$$0<a_{m,k} < {1\over k\log m} \to 0 \quad {\rm as} \quad m\to \infty \quad {\rm for} 
\quad k=1,2,\ldots.$$
Thus, the sufficient conditions are satisfied in Toeplitz' theorem, and we conclude that the limit in (1.3) holds in the particular choice when 
$t=m\in \N$. 

(ii) Second, given any real number $t>3$, let $m:= [t]$, the integer part of $t$. We use (4.3) and the Second Mean-Value Theorem again to obtain 
$$\tau(t) \log t - \tau(m) \log m = \int^t_m {s(u)\over u} du\leqno(4.5)$$
$$={1\over m} \int^\xi_m s(u) du - {1\over t} \int^t_\xi s(u) du$$
$$=\Big({1\over m} - {1\over t}\Big) \xi\sigma(\xi) - \sigma(m) + \sigma(t), \quad m<\xi<t.$$
By (1.1) and (4.5), we get 
$$|\tau(t) \log t - \tau(m) \log m| = {t-m\over mt} \xi|\sigma(\xi)| + |\sigma(t) - \sigma(m)|$$
$$\le {1\over m}|\sigma(\xi)| + |\sigma(t) - \sigma(m)| \to 0\quad {\rm as} \quad m\to \infty.$$
Hence we conclude that 
$$\lim_{t\to \infty} \tau(t) = \lim_{t\to \infty} \tau(m) {\log m\over \log t} = A, \quad {\rm where}\quad m:= [t].$$

(iii) Third, to see that the converse implication is not true in general, we consider the function $s$ 
defined by 
$$s(t):=\cases{me^{2^m} \quad &if\quad $t\in [e^{2^m}, e^{2^m}+1], m=1,2,\ldots$;\cr
0 &otherwise on\quad $[1, \infty)$.\cr}$$

We claim that this function $s$ cannot be summable $(C, 1)$ to any finite number $A$. To this effect, we recall that 
if we had (1.1), then for any number $a>0$ we would have 
$${1\over t} \int^{t+a}_t s(u) du = {t+a\over t} \sigma(t+a) - \sigma(t) \to 0 \quad 
{\rm as} \quad t\to \infty.$$
But for $t:= e^{2^m}$ and $a:= 1$, we have 
$${1\over t} \int^{t+1}_t s(u) du = e^{-2^m} \int^{t+1}_t me^{2^m} du = m\not\to 0\quad {\rm as} \quad m\to \infty.$$
Consequently, for this function $s$ the limit (1.1) cannot exist with any finite number $A$. 

On the other hand, if $t$ is such that 
$$e^{2^{m-1}} \le t<e^{2^m}, \quad m=1,2,\ldots;$$
then we estimate as follows: 
$$\eqalign{0\le \tau(t) &\le {1\over 2^{m-1}} \sum^{m-1}_{k=1} \int^{e^{2^k}}_{e^{2^{k-1}}} {s(u)\over u} du\cr
&\le {1\over 2^{m-1}} \sum^{m-1}_{k=1} k\to 0 \quad {\rm as} \quad 
t \to \infty.\cr}$$
This proves that the limit in (1.3) exists with $A=0$. 

The proof of Theorem 3 is complete. 

We note that summability $(L, 2)$ is more effective than summability $(L, 1)$. This can be proved 
in an analogous way as Theorem 3 was proved above. 
We refer to [5, on p. 382], where an analogous result is proved for the logarithmic 
mean $\tau_2(n)$ of second order of a numerical sequence $(s_k)$ (see also (5.2) below). 

\vglue1cm

\centerline{5. {\it Summability} $(L, 1)$ {\it of numerical sequences} }

The above methods of summability are the nondisrete ones of the methods of logarithmic summability of 
numerical sequences $(s_k) = (s_k : k=1,2,\ldots)$ of complex numbers. We recall that a sequence 
$(s_k)$ is said to be logarithmic summable of order 1 (see in [5], where the term `harmonic summable of order 1' was used), 
briefly: summable $(L, 1)$, if there exists some $A\in \C$ such that 
$$\lim_{n\to \infty} {1\over \ell_n} \sum^n_{k=1} {s_k\over k} = A, \quad {\rm where}\quad 
\ell_n := \sum^n_{k=1} {1\over k} \sim \log n, \leqno(5.1)$$
where for two 
sequences $(a_n)$ and $(b_n)$ of positive numbers we write $a_n \sim b_n$ if 
$$\lim_{n\to \infty} {a_n\over b_n} = 1.$$ 

We note that the sequence $(s_k)$ is said to be logarithmic summable of order 2 (see also in [5]), briefly: 
summable $(L, 2)$, if there exists some $A\in \C$ such that 
$$\lim_{n\to \infty} \tau_2 (n) := {1\over \ell_n (2)} \sum^n_{k=1} {s_k \over k\ell_k},
\quad {\rm where}\quad 
\ell_n (2) := \sum^n_{k=1} {1\over k\ell_k} \sim \log \log n.
\leqno(5.2)$$

It is clear that if the ordinary limit 
$$\lim_{n\to \infty} s_n = A\leqno(5.3)$$
exists, then the limit in (5.1) also exists with the same $A$. Even more is true (see, e.g., in [5, on p. 376]): 
If a sequence $(s_k)$ is such that the finite limit 
$$\lim_{n\to \infty} {1\over n} \sum^n_{k=1} s_k = A$$
exists, then the limit in (5.1) also exists with the same $A$. The converse implication is not true in general. 

We note that if the finite limit in (5.1) exists, then the limit in (5.2) also exists with the same $A$ 
(see also in [5, on p. 382]. Again, the converse implication is not true in general. 

Now, the discrete analogue of Theorem 1 reads as follows. 

\bigskip

THEOREM 4. {\it If a sequence $(s_k)$ of real numbers is summable $(L, 1)$ to some $A\in \R$, then the 
ordinary limit} (5.3) {\it exists if and only if 
$$\limsup_{\lambda\to 1+} \ \liminf_{n\to \infty} {1\over ([n^\lambda]-n) \ell_n} \ 
\sum^{[n^\lambda]}_{k=n+1} {s_k-s_n\over k} \ge 0\leqno(5.4)$$
and
$$\limsup_{\lambda\to 1-} \ \liminf_{n\to \infty} \ 
{1\over (n-[n^\lambda])\ell_n} \sum^n_{k=[n^\lambda]+1} {s_n - s_k\over k} \ge 0, \leqno(5.5)$$
where by {\rm [ $\cdot$ ]} we denote the integer part of a real number, and $\ell_n$ is defined in} (5.1). 

Analogously to (2.3), we say that a sequence $(s_k)$ of real numbers is slowly decreasing with respect to 
summability $(L, 1)$ if for every $\e>0$ there exist a natural number $n_0 = n_0 (\e)$ and a real number 
$\lambda=\lambda(e) > 1$ such that 
$$s_k - s_n \ge -\e\quad {\rm whenever}\quad n_0 \le n < k \le n^\lambda.\leqno(5.6)$$

It is easy to check (cf. (2.4)) that a sequence $(s_k)$ is slowly decreasing with respect to summability $(L, 1)$ if and 
only if 
$$\lim_{\lambda\to 1+} \ \liminf_{n\to \infty} \ \min_{n<k\le n^\lambda} (s_k - s_n) \ge 0. \leqno(5.7)$$

Clearly, if a sequence $(s_k)$ is slowly decreasing with respect to summability $(L, 1)$, then both 
conditions (5.4) and (5.5) are satisfied. Thus, the next corollary is an immediate consequence of Theorem 4. 

\bigskip

COROLLARY 3. {\it Suppose a sequence $(s_k)$ of real numbers is slowly decreasing with respect to summability $(L, 1)$. 
If $(s_k)$ is summable $(L, 1)$ to some $A\in \R$, then the ordinary limit} (5.3) {\it also exists.} 

This corollary was earlier proved   by Kwee [3, Lemma 3] in a different way. We note that the definition of slow decrease of a sequence 
$(s_k)$ is formally different in [3] from the definitions given in (5.6) and (5.7) above. 

\bigskip

REMARK 2. According to Kwee's definition  in  [3,  see it as a condition in both Theorem A and Lemma 3], a sequence of real 
numbers $(s_k)$ is said to be 
slowly decreasing if
$$\liminf(s_k - s_n)\ge 0\quad {\rm whenever}\quad k>n\to \infty\quad {\rm and} \quad 
{\log k\over \log n} \to 1. \leqno(5.8)$$

We claim that definition (5.8) is equivalent to the one in (5.6) (as well as to the one in (5.7)). The implication (5.8) $\Rightarrow$ 
(5.6) is trivial. To justify the converse implication (5.8) $\Rightarrow$ (5.6), let $\lambda>1$ be 
arbitrarily close to 1 and set $\e:= \log \lambda$. By (5.8), there exists $n_0 = n_0(\lambda) > 1$ such that 
$$s_k-s_n\ge -\e \quad {\rm whenever} \quad k>n\ge n_0\quad {\rm and}\quad 
0<\log {\log k\over \log n} \le \log \lambda = \e.$$
Now, the equivalence of the two definitions claimed above is obvious. 

Next, the discrete analogue of Theorem 2 reads as follows. 

\bigskip

THEOREM 5. {\it If a sequence $(s_n)$ of complex numbers is summable $(L, 1)$ 
to some $A\in \C$, then the ordinary limit} (5.3) {\it exists if and only if} 
$$\lim_{\lambda\to 1+} \ \limsup_{t\to \infty} \ \Big|{1\over ([n^\lambda]-n) \ell_n} 
\sum^{[n^\lambda]}_{k=n+1} {s_k-s_n\over k}\Big| = 0.\leqno(5.9)$$

Analogously to (2.6), we say that a sequence $(s_k)$ of complex numbers is slowly oscillating with respect to 
summability $(L, 1)$  if for every $\e>0$ there exist 
$n_0 = n_0 (\e) > 1$ and $\lambda=\lambda(\e) > 1$ such that 
$$|s_k - s_n| \le \e \quad {\rm whenever} \quad n_0 \le n < k \le n^\lambda.\leqno(5.10)$$

It is easy to check that a sequence $(s_k)$ is slowly oscillating with respect to summability 
$(L, 1)$ if and only if 
$$\lim_{\lambda\to 1+} \ \limsup_{n\to \infty} \ \max_{n<k\le n^\lambda} |s_k-s_n|=0.\leqno(5.11)$$

\bigskip

REMARK 3. The concept of slow oscillation with respect to summability $(L,1)$ is not defined in [3]. However, analogously 
to (5.8), a sequence $(s_k)$ of complex numbers may be called to be slowly 
oscillating if 
$$\lim(s_k-s_n) = 0\quad {\rm whenever} \quad k> n\to \infty\quad {\rm and}\quad 
{\log k\over \log n} \to 1.\leqno(5.12)$$
A reasoning similar to the one in Remark 2 gives that the definitions (5.12) and (5.10) (as well as (5.11)) are equivalent. 

It is clear that if a sequence $(s_k)$ is slowly oscillating with respect to $(L,1)$, then condition 
(5.11) is satisfied. Thus, the next corollary is an immediate consequence of Theorem 5. 

COROLLARY 4. {\it Suppose a sequence $(s_k)$ of complex numbers is slowly oscillating with respect to summability 
$(L,1)$. If $(s_k)$ is summable $(L,1)$ to some $A\in \C$, then the ordinary limit} (5.3) {\it also exists}. 

The proofs of Theorems 4 and 5 run along similar lines to those of Theorems 1 and 2, respectively; 
while the key ingredient is provided by the following

\bigskip

LEMMA 2. (i) {\it For all $\lambda>1$ and large enough $n$, that is when $[n^\lambda]>n$, we have the representation} 
$$s_n - \tau_n = {\ell_{[n^\lambda]}\over \ell_{[n^\lambda]}-\ell_n} (\tau_{[n^\lambda]} - \tau_n) - 
{1\over \ell_{[n^\lambda]}-\ell_n} \sum^{[n^\lambda]}_{k=n+1} {s_k-s_n\over k}.$$

(ii) {\it For all $0<\lambda<1$ and large enough $n$, that is when 
$n> [n^\lambda]$, we have} 
$$s_n - \tau_n = {\ell_{[n^\lambda]}\over \ell_n - \ell_{[n^\lambda]}} (\tau_n - \tau_{[n^\lambda]}) + 
{1\over \ell_n - \ell_{[n^\lambda]}} \sum^n_{k=[n^\lambda]+1} {s_n -s_k\over k}.$$

{\it Proof}. Performing steps analogous to those in the proof of 
Lemma 1 yields the above representations. 

\vfill\eject

\centerline{REFERENCES}

\item{[1]} G. H. HARDY, Theorems relating to the summability and 
convergence of slowly oscillating sequences, {\it Proc. London Math. Soc.} (2), {\bf 8}(1910), 301-320. 

\item{[2]} G. H. HARDY, {\it Divergent Series}, Clarendon Press, Oxford, 1949. 
\item{[3]} B. KWEE, A Tauberian theorem for the logarithmic method of summation, 
{\it Proc. Camb. Phil. Soc.}, {\bf 63} (1967), 401-405. 
\item{[4]} E. L. LANDAU, \"Uber die Bedeutung einiger neuen Grenzwerts\"atze der Herren Hardy and Axel, 
{\it Prace Mat.-Fiz. Analysis} {\bf 21} (1910), 97-177. 
\item{[5]} F. M\'ORICZ, On the harmonic averages of numerical sequences, 
{\it Archiv Math. (Basel)} {\bf 86}(2006), 375-384. 
\item{[6]} R. SCHMIDT, \"Uber divergente Folgen und Mittelbildungen, {\it Math. Z.} {\bf 22}(1925), 89-152. 
\item{[7]} E. C. TITCHMARSH, {\it Introduction to the Theory of Fourier Integrals}, Clarendon Press, 
Oxford, 1937. 
\item{[8]} A. ZYGMUND, {\it Trigonometric Series}, Vol. I., Cambridge Univ. Press, 1959. 

\bye